\UseRawInputEncoding
\documentclass{amsart}
\usepackage{amsmath}
\usepackage{amssymb}
\newtheorem{thm}{Theorem}
\newtheorem{athm}{Theorem}[section]
\newtheorem{cor}[thm]{Corollary}
\newtheorem{lem}[thm]{Lemma}
\newtheorem{alem}[athm]{Lemma}
\newtheorem{prop}[thm]{Proposition}
\newtheorem{cons}[thm]{Construction}
\newtheorem{ex}[thm]{Example}
\theoremstyle{definition}
\newtheorem{defn}[thm]{Definition}
\newtheorem{adefn}[athm]{Definition}
\newtheorem{notn}[thm]{Notation}
\theoremstyle{remark}
\newtheorem{rem}[thm]{Remark}
\newtheorem{arem}[athm]{Remark}
\long\def\Thm#1{\begin{thm} #1 \end{thm}}
\long\def\aThm#1{\begin{athm} #1 \end{athm}}
\long\def\Cor#1{\begin{cor} #1 \end{cor}}
\long\def\Lem#1{\begin{lem} #1 \end{lem}}
\long\def\aLem#1{\begin{alem} #1 \end{alem}}

\long\def\aDef#1{\begin{adefn} #1 \end{adefn}}

\long\def\Rem#1{\begin{rem} #1 \end{rem}}
\long\def\aRem#1{\begin{arem} #1 \end{arem}}
\long\def\Ex#1{\begin{ex} #1 \end{ex}}

\def\Sect{\section}

\long\def\Ref#1#2#3#4#5#6{
\bibitem{#1}
{\rm #2,}
\textit{#3.}
{\rm #4}
\textbf{#5}
{\rm #6.}
}
\long\def\Refb#1#2#3#4{
\bibitem{#1}
{\rm #2,}
\textit{#3.}
#4.
}
\def\Rr{{\mathbb R}}
\def\Ff{{\mathbb F}}

\def\into{\hookrightarrow}
\def\leq{\leqslant}
\def\geq{\geqslant}

\def\st{\mid}

\def\phi{\varphi}

\def\cl#1{\overline{#1}}

\def\Ff{\mathbb{F}}

\begin{document}
\title{A hyperplane Ham Sandwich theorem}
\author{M.~C.~Crabb}
\address{%
Institute of Mathematics\\
University of Aberdeen \\
Aberdeen AB24 3UE \\
UK}
\email{m.crabb@abdn.ac.uk}
\begin{abstract}
We give a direct proof of a result due to Karasev (2008),
Karasev-Matschke (2014) and Schnider-Sober\'on (2023).
Given $m+1$ Borel probability measures on the space of affine
hyperplanes in a real vector space $V$ of dimension $m+1$,
there exist a line $L$ through the origin in $V$ and a point
$v\in L$ such that at least half of the hyperplanes, as counted
by any of the measures, meet or are parallel to each of the two
closed rays in $L$ meeting at $v$.
\end{abstract}
\date{January 2024, revised March 2024}
\maketitle

Let $V$ be a Euclidean vector space of dimension $m+1$,
where $m\geq 0$.
The affine hyperplanes in $V$ are parametrized by equivalence
classes $[f ,y]$ of pairs
$(f ,y)$, where $f\in V^*-\{ 0\}$ is a non-zero linear form
on $V$ and $y\in\Rr$,
and $[f, y]=[tf ,ty]$ for $t\in\Rr^\times$, the class $[f,y]$ giving
the hyperplane $\{ v\in V\st f(v)=y\}$.
Write $H^*$ for this space; it is the total space of the dual
Hopf bundle over the real projective space $P(V^*)$ of $V^*$.
Let us write $H$ for the total space of the Hopf line
bundle over $P(V)$ with points written as $[e,x]$, where 
$e\in V-\{0\}$ and $x\in\Rr$, $[e,x]=[te,x/t]$ for $t\in\Rr^\times$,
so that $L=\Rr e$ is a line in $V$ and $v=xe\in L$
is a point on the line.
(One can think of $H$ as the space obtained from $V$
by `blowing up' the point $0$ in the language of Algebraic Geometry.)

In this note we present a direct proof of
the following ``hyperplane Ham Sandwich theorem'' due
in different forms to Karasev \cite{karasev},
Karasev-Matschke \cite{matschke} and Schnider-Sober\'on \cite{soberon}.
\Thm{\label{main}
Suppose that $\mu_0,\ldots ,\mu_l$, $l\geq 0$,
are Borel probability measures
on the space $H^*$ of affine hyperplanes in the $(m+1)$-dimensional
real vector space $V$.

If $l \leq m$, then there is a point $[e,x]\in H$ in the total 
space of
the Hopf bundle over $P(V)$ with the property $(*)$ that
$$\textstyle
\mu_j\{ [f ,y] \in H^*\st y f (e)\geq x f(e)^2\}
\geq \frac{1}{2}
\text{\quad and\quad}
\mu_j\{ [f ,y]\in H^* \st y f (e)\leq x f(e)^2\}
\geq \frac{1}{2}
$$
for $j=0,\ldots ,l$.
}
The case in which the measures are absolutely continuous  
with compact supports is contained in \cite[Theorem 3]{karasev};
a compactified (projective) version for arbitrary measures,
equivalent to Theorem \ref{main},
is included in
\cite[Corollary 5.1]{matschke}; and the following specialization
to the discrete case is in \cite[Corollary 2]{soberon}.
\Cor{\label{schnider}
Suppose that $M_0,\ldots ,M_m$ are non-empty finite sets of
affine hyperplanes in $V$ in the non-zero real vector space 
$V$ of dimension $m+1$.

Then there exist a line $L\subseteq V$ through zero and a point 
$v\in L$ such that, for each of the two closed rays $R$
in $L$ starting
at $v$ and each $j=0,\ldots ,m$, the number of hyperplanes
in $M_j$ that meet $R$ or are parallel to $L$ is greater
than or equal to $\# M_j/2$.
}
\begin{proof}
Let $\mu_j$ in Theorem \ref{main} be the discrete measure
taking the value $1/\# M_j$ on each of the points of $M_j$.
Then $L=\Rr e$ and $v=xe$.
\end{proof}
\Ex{Let $e_0,\ldots ,e_m$ be a basis of $V$ and $f_0,\ldots ,f_m$ the dual basis of $V^*$. Consider the $1$-element subsets
$M_j=\{ [f_j,1]\}$ in Corollary \ref{schnider}. 
Then the possible solutions
$(L,v)$ are parametrized by the $2^{m+1}-1$ non-empty
subsets $J\subseteq \{ 0,\ldots ,m\}$ with
$v=(\sum_{j\in J}e_j)/\# J$.
{\rm Since the hyperplane $[f_0+\ldots +f_m,0]$ is not parallel
to any solution line $L$ and does not pass through any solution
point $v$, we see that the restriction $l\leq m$ in Theorem
\ref{main} is optimal.}
}
The proof is, in fact, only slightly harder than the proof of the classical Ham Sandwich theorem, which we include for comparison in an Appendix. It will depend upon the following input from 
Algebraic Topology, which is a consequence of the fact that
the power $e(H)^l$ of
the $\Ff_2$-cohomology Euler class $e(H)\in H^1(P(V);\,\Ff_2)$
of the Hopf line bundle $H$ over the real projective space $P(V)$
of $V$ is zero if and only if $l > m$.
\Lem{\label{euler}
Let $V$ be a Euclidean vector space of dimension $m+1$.
Suppose that $s$ is a nowhere zero section of
the pullback of $\Rr^{l+1}\otimes H$ to the unit disc
bundle $D(H)$ over $P(V)$ with the property that
the restriction of $s$ to the sphere bundle $S(H)$
is homotopic, through nowhere zero sections,
to the diagonal inclusion $S(H)\into \Rr^{l+1}\otimes H$.
Then $m <l$.
\qed
}
\begin{proof}[Proof of Theorem \ref{main}]
Choose a Euclidean inner product on $V$
and write $S(V)$ for the unit sphere in $V$.
This allows us now to take representative points of $H$ and $H^*$ as
$[e,x]$ and $[f,y]$ with $e\in S(V)$ and $f\in S(V^*)$
of unit norm.

We distinguish two cases:

\par\noindent (i).
There exists some $\rho >0$ such that for each $e\in S(V)$
$$\textstyle
\mu_j\{ [f ,y]\in H^*\st |f(e)|\leq |y|/\rho \}
<\frac{1}{2}
$$
for some $j$.

\par\noindent (ii).
Property (i) does not hold, and so there exists
a sequence $e_n$, $n\geq 1$, in $S(V)$ such that
$$\textstyle
\mu_j\{ [f ,y]\in H^*\st |f (e_n)|\leq |y|/n \}
\geq \frac{1}{2}
$$
for all $j=0,\ldots ,l$.

In case (i) we shall show that either (a)
there exists
some $(e,x)$ with $|x|<\rho$ such that property $(*)$
holds for all $j=0,\ldots ,l$, or (b) $l >m$.
Suppose, then, that (i) holds but there is no
$(e,x)$ with $|x|<\rho$ satsifying $(*)$.

Define $\Omega_j$, for $j=0,\ldots ,l$, to be the subset
$$\textstyle
\{ [e,x]\in H\st 
\mu_j\{ [f ,y] \in H^*\st y f (e)\geq x f(e)^2\}
<\frac{1}{2}
$$
$$\textstyle\qquad
\text{\ or\quad }
\mu_j\{ [f ,y]\in H^* \st y f(e)\leq x f (e)^2\}
<\frac{1}{2}\}
$$
of $H$,
and define a section $s_j$ of the pullback of $H$ to
$\Omega_j$ by 
$$
s_j[e,x]=
\begin{cases}
[e,1]&\text{if $\mu_j\{ [f ,y]\in H^* \st y f (e)\geq x f (e)^2\}
<\frac{1}{2}$},\\
[e,-1]&\text{if $\mu_j\{ [f ,y]\in H^* \st y f (e)\leq x f (e)^2\}
<\frac{1}{2}$}.
\end{cases}
$$
(The two cases are mutually exclusive
because the sum of the two measures
must be greater than or equal to $1$.)

If $(e,x)\in S(V)\times\Rr$ satisfies 
$\mu_j\{ [f ,y]\in H^* \st y f (e)\geq x f (e)^2\}
<\frac{1}{2}$, then a routine argument
(given as Lemma \ref{detail} below)
shows that there is an open neigbourhood
$N$ of $(e,x)$ such that 
$\mu_j\{ [f ,y]\in H^* \st y f (e')\geq x' f (e')^2\}
<\frac{1}{2}$ for all $(e',x')\in N$.
Hence the subspaces $\Omega_j\subseteq H$ are open
and the sections $s_j$ are continuous.

If $y f (e)\geq \rho f (e)^2$, then
$\rho f(e)^2 \leq |y|\, |f (e)|$ and so
$|f(e)| \leq |y|/\rho$.
Thus $[e,\rho]\in\Omega_j$ and $s_j[e,\rho]=[e,1]$.

Write $\Omega_j^\rho =(\Omega_j\cap D_\rho (H))-
\{[e,\rho]\in\Omega_j\st s_j[e,\rho]=[e,-1]\}$.
By assumption, the closed disc bundle
$D_\rho (H)$ of radius $\rho$ is contained in 
$\bigcup_{j=0}^l \Omega_j$. So $(\Omega_j^\rho)_{j=0}^l$ 
is an open cover
of $D_\rho (H)$ and we can choose a subordinate
partition of unity $(\phi_j)_{j=0}^l$.
Now the section $s=(\phi_js_j)$ of $\Rr^{l+1}\otimes H$
over $D_\rho (H)$ is nowhere zero and its restriction
to $S_\rho (H)$ is linearly homotopic through nowhere
zero sections to the diagonal inclusion
$S_\rho(H)\into \Rr^{l+1}\otimes H$.
It follows from Lemma \ref{euler} that $m < l$.

This deals with case (i). Suppose now that (ii) holds.
By the compactness of the sphere, we may assume that the
sequence $e_n$ converges to some $e\in S(V)$. Then
$$\textstyle
\mu_j\{ [f ,y]\in H^*\st f (e)=0 \}
\geq \frac{1}{2}
$$
for all $j=0,\ldots ,l$, because
$$
\{ [f ,y]\in H^*\st f (e)=0 \}=\bigcap_{n\geq 1}
\{ [f ,y]\in H^*\st |f (e_n)|\leq |y|/n \}.
$$
So $[e,0]$, and, indeed, $[e,x]\in H$ for any $x\in\Rr$,
satisfies the condition $(*)$ for $j=0,\ldots ,l$.
\end{proof}
\Lem{\label{detail}
Suppose that $(e,x)\in S(V)\times\Rr$ in the proof of Theorem 
\ref{main}
satisfies 
$\mu_j\{ [f ,y]\in H^* \st y f (e)\geq x f (e)^2\}
=\frac{1}{2}-3\delta$, where $\delta >0$. 
Then there is an open neigbourhood
$N$ of $(e,x)$ such that 
$\mu_j\{ [f ,y]\in H^* \st y f (e')\geq x' f (e')^2\}
<\frac{1}{2}-\delta$ for all $(e',x')\in N$.
}
\begin{proof}
There exist some $\epsilon >0$ such that
$\mu_j\{ [f ,y]\in H^* \st y f(e)> x f (e)^2-\epsilon\}
<\frac{1}{2}-2\delta$ and
some $R>0$ such that
$\mu_j\{ [f ,y]\st \| f\| =1,\, |y|\geq R\} <\delta$.
So there is a neighbourhood $N$ of $(e,x)$ in $S(V)\times\Rr$ 
such that,
if $y f (e')\geq x' f (e')^2$, $[e',x']\in N$
and $\|f\| =1$, $|y|<R$, then $y f (e)>x f (e)^2-\epsilon$. 
And then $\mu_j\{ [f ,y]\in H^* \st y f (e')\geq x' f (e')^2\}
< \frac{1}{2}-2\delta +\delta=\frac{1}{2}-\delta $.
\end{proof}
\Rem{For an example of case (ii),
take $l=1$, $m\geq 1$, and,
for $f\in S(V^*)$, let $\mu_0$ and $\mu_1$ be the
discrete measures concentrated at the single points
$[f,0]$ and $[f,1]$ respectively as in Corollary \ref{schnider}.
The solutions $(L,v)$ consist of the lines $L$ in the kernel
of $f$ and arbitrary $v\in L$. The line $L$ meets 
the hyperplane $[f,0]$ at $v$ if and only if $v=0$,
and the hyperplane $[f,1]$ if and only $f(v)=1$,
but is parallel to both hyperplanes.
}
\begin{appendix}
\Sect{The Ham Sandwich theorem}
Here is the analogous proof of the classical Ham Sandwich theorem.
A much more general result, established by similar methods
and to which this appendix may serve as an introduction,
can be found in \cite[Proposition 9.2 and Theorem 9.3]{BC}.
\aThm{Suppose that $\nu_0,\ldots ,\nu_l$ are Borel probability
measures on the real vector space $V$ of dimension $m+1$. 
If $l \leq m$, then there exists an affine hyperplane
$[f ,y]\in H^*$ such that 
$$\textstyle
\nu_j\{ v\in V\st y\geq f (v)\}
\geq \frac{1}{2}
\text{\quad and\quad}
\nu_j\{ v\in V\st y\leq f (v)\}
\geq \frac{1}{2}
\leqno{(**)}
$$
for $j=0,\ldots ,l$.
}
\begin{proof}
Fix a Euclidean inner product on $V$.
Choose $\rho >0$ such that
$$\textstyle
\nu_j\{ v\in V\st \| v\| \geq \rho\} <\frac{1}{2}
$$
for $j=0,\ldots ,l$.
Suppose that there is no $[f ,y]\in H^*$ with
$\| f\| =1$ and $|y| <\rho$ satisfying $(**)$.
We shall show that $l>m$.

Define $\Omega_j$, for $j=0,\ldots ,l$, to be the subset
$$\textstyle
\{ [f ,y]\in H^*\st \nu_j\{ v\in V\st f (v)\geq y\}
<\frac{1}{2}
\text{\quad or\quad }
\nu_j\{ v\in V\st f (v)\leq y\}
<\frac{1}{2}\}
$$
of $H^*$.
A section $s_j$ of the pullback of $H^*$ to $\Omega_j$ is
defined, for $f\in S(V^*)$, by
$$
s_j[f ,y]=\begin{cases}
[f,1]&\text{if $\nu_j\{ v\in V\st f (v)\geq y\}
< \frac{1}{2}$,}\\
[f,-1]&\text{if $\nu_j\{ v\in V\st f (v)\leq y\}
< \frac{1}{2}$.}
\end{cases}
$$
As in the proof of Theorem \ref{main} we see that $\Omega_j$
is open and $s_j$ is continuous.

If $f (v)\geq y$, then $\| v\| \geq y$. So $[f ,\rho ]\in
\Omega_j$ for all $j$ and $s_j[y,\rho]=[y,1]$.

We are assuming that $D_\rho (H^*)\subseteq \bigcup_j\Omega_j$.
So $(\Omega_j\cap D_\rho (H^*))$ is an open cover of $D_\rho (H^*)$.
Choose a partition of unity $(\phi_j)$ subordinate to 
this cover.
Then $s=(\phi_js_j)$ is a nowhere zero section of
$\Rr^{l+1}\otimes H^*$ over $D_\rho (H^*)$,
and $s$ restricts on $S_\rho(H^*)$ to a section that
is linearly homotopic through nowhere zero sections
to the diagonal inclusion. So, by Lemma \ref{euler}, 
we conclude that $m <l$.
\end{proof}
\aRem{Suppose that $\rho>0$ in the proof above is required only to
satisfy the property that, for each $f\in S(V^*)$, there is some
$j$ such that
$$\textstyle
\nu_j\{ v\in V\st f(v)\geq \rho\} <\frac{1}{2}.
$$
Then the argument can be adapted, 
as in the proof of Theorem \ref{main},
to show that, if $l\leq m$, there is a solution
$[f,y]$ with $\| f\| =1$ and $|y| <\rho$.
}
\Sect{A parametrized version}
The earlier results can be extended, along the lines of similar
generalizations in \cite{tmna, BC},
to a fibrewise version
in which the Euclidean space $V$ is replaced by
a real vector bundle $E$ of dimension $m+1$ over a compact
ENR (Euclidean Neighbourhood Retract) $B$. 
We may assume that
$E$ is equipped with a Euclidean inner product
and consider the real projective bundle
$P(E)\to B$ and the Hopf line bundle $H$ over $P(E)$.
At a point $b\in B$, the fibre is the real projective
space $P(E_b)$ of the fibre $E_b$ of $E$; we write
$H_b$ for the restriction of $H$ to $P(E_b)$.
\aDef{A family $(\mu^b)_{b\in B}$ of Borel probability
measures $\mu^b$ on the space $H^*_b$ of affine hyperplanes
in $E_b$, parametrized by the compact ENR $B$,
is said to be {\it continuous}
if, for each continuous function $\phi : H^*\to \Rr$ with
compact support, the function $b\mapsto \int \phi_b\, {\rm d}\mu^b :
B\to\Rr$ which integrates the restriction $\phi_b$ of $\phi$
to $H_b^*$ with respect to the measure $\mu^b$ is continuous.
}
\aRem{\label{continuity}
Suppose that $E$ is trivial: $E=B\times V$ and that the family 
$(\mu^b)$ is continuous. The fibre $H^*_b$ of $H^*$ at $b$
is thus independent of $b\in B$; let us write it as $F$
and regard $(\mu^b)$ as a family of measures in $F$.
Let $W\subseteq F$ be open
and $t\in\Rr$. Then $\{ b\in B \st \mu^b(W)>t\}$ is open
in $B$.
}
\begin{proof}
Suppose that $\mu^b(W)>t$.
Then there is a continuous function $\chi : F=H^*_b \to [0,1]$ 
with compact support
$\cl{\chi^{-1}(0,1]}\subseteq W$ such that
$\mu^b(W)\geq  \int \chi {\rm d}\mu^b>t$.
By the continuity of the family, applied to the composition
$\phi : H^*=B\times F \to \Rr$ of $\chi$ and the projection,
$\mu^{b'}(W)\geq  \int \chi {\rm d}\mu^{b'}>t$
for all $b'$ in some neighbourhood of $b$.
\end{proof}
The generalization of Theorem \ref{main} involves the Stiefel-Whitney
classes $w_j(-E)\in H^j(B;\,\Ff_2)$ of the stable complement 
$-E$ of $E$.
\aThm{\label{fw_main}
Let $E\to B$ be a Euclidean vector bundle
of dimension $m+1$ over a compact ENR $B$.
Suppose that $\mu_j^b$ $(b\in B)$ $j=0,\ldots ,l$, 
are continuous families of Borel probability measures on the
spaces $H^*_b$ of affine hyperplanes in $E_b$.
If $l\geq m$ and $w_{l-m}(-E)\not=0$, then there is a point $b\in B$
of the base and $[e,x]\in H_b$ in the total space of
the Hopf bundle over $P(E_b)$ with the property $(*)$ that
$$\textstyle
\mu_j^b\{ [f ,y] \in H_b^*\st y f (e)\geq x f(e)^2\}
\geq \frac{1}{2}
\text{\quad and\quad}
\mu_j^b\{ [f ,y]\in H_b^* \st y f (e)\leq x f(e)^2\}
\geq \frac{1}{2}
$$
for $j=0,\ldots ,l$.
}
The proof requires the following criterion for the vanishing of
an $\Ff_2$-cohomology Euler class.
\aLem{\label{fw_euler}
Suppose that $s$ is a nowhere zero section of
the pullback of $\Rr^{l+1}\otimes H$ to the unit disc bundle
$D(H)$ over $P(E)$ such that the restriction of $s$ to
$S(H)\subseteq D(H)$ is homotopic through nowhere zero sections
to the diagonal inclusion $S(H)\into \Rr^{l+1}\otimes H$.
Then $l>m$ and each Stiefel-Whitney class $w_j(-E)\in H^j(B;\,\Ff_2)$
for $j\geq l-m$ is zero.
}
\begin{proof}
We recall from, for example, \cite[Proof of Proposition 4.1]{jan} and
\cite[Proof of Proposition 2.4]{BC} that
$$
H^*(P(E);\,\Ff_2)=H^*(B;\,\Ff_2)[T]/(w_{m+1}(E)+w_m(E)T+\cdots
+T^{m+1}),
$$ 
and $e(H)^l$ is represented by $d_0+d_1T+\cdots +d_mT^m$, where 
$d_j=w_{l-j}(-E)+w_1(E)w_{l-j-1}(-E)+\cdots + w_{m-j}(E)w_{l-m}(-E)$.

The existence of the section $s$ implies that
the Euler class $e(H)^l\in H^l(P(E);\, \Ff_2)$ of $\Rr^l\otimes H$
is zero. So $l>m$ and $w_{l-m}(-E)=0$,
$w_{l-m+1}(-E)=0$, $\ldots$, $w_l(-E)=0$.
For any $l'\geq l$, we have $e(H)^{l'}=0$, too, and thus
$w_{l'}(-E)=0$.
\end{proof}
\begin{proof}[Outline proof of Theorem \ref{fw_main}]
Following closely the proof of Theorem \ref{main}, we distinguish two cases:

\par\noindent 
(i). There exists a real number
$\rho >0$ such that for each $b\in B$ and $e\in S(E_b)$ there is some $j$ such that
$$\textstyle
\mu_j^b\{ [f,y]\in H^*_b \st | f(e)|\leq |y|/\rho\} <\frac{1}{2}.
$$
\par\noindent 
(ii). Property (i) fails to hold.

Assume first that (i) holds, but there is no point $[e,x]\in H_b$
having the property $(*)$. If we can construct a section $s$ as in
Lemma \ref{fw_euler}, the result will follow.

For $b\in B$,
define $\Omega_j^b$, $j=0,\ldots ,l$, to be the subset
$$\textstyle
\{ [e,x]\in H_b\st 
\mu_j^b\{ [f ,y] \in (H_b)^*\st y f (e)\geq x f(e)^2\}
<\frac{1}{2}
$$
$$\textstyle\qquad
\text{\ or\quad }
\mu_j^b\{ [f ,y]\in (H_b)^* \st y f(e)\leq x f (e)^2\}
<\frac{1}{2}\}
$$
of $H_b$. Then put 
$\Omega_j =\bigcup_{b\in B} \Omega_j^b\subseteq H$.

If $(e,x)\in S(E_b)\times\Rr$ satisfies
$\mu^b_j\{ [f,y]\in H^*_b\st yf(e)\geq xf(e)^2\} <\frac{1}{2}$,
then one can show, by adapting the proof of Lemma \ref{detail}
and using Remark \ref{continuity} 
(given that $E$ is trivial in a neighbourhood of $b$), that
there is an open neighbourhood $N$ of $(b,(e,x))$ in $H^*$
such that $\mu^{b'}_j\{ [f,y]\in H^*_{b'}\st yf(e')\geq x'f(e')^2\} <\frac{1}{2}$ for all $(b',(e',x'))\in N$.
Continuous sections $s_j$ of $H$ over $\Omega_j$
and the required section $s$ of $\Rr^{l+1}\otimes H$ can then be
constructed by following the proof of Theorem \ref{main}.

In case (ii), since $S(E)$ is compact, there is a convergent
sequence $e_n\in S(E_{b_n})$, $n\geq 1$, such that
$$\textstyle
\mu_j^{b_n}\{ [f,y]\in H_{b_n}^* \st |f(e_n)| \leq |y|/n\} 
\geq\frac{1}{2}
$$
for $j=0,\ldots ,l$.
As the bundle $E$ is trivial in a neighbourhood of 
$b=\lim_{n\to\infty} b_n$,
there is no loss of generality in now supposing that  $E$ is trivial,
say $E=B\times V$, where $V=E_b$. 
Let $e=\lim_{n\to\infty} e_n\in S(V)$.
The proof will be completed by showing that
$$\textstyle
\mu_j^b\{ [f,y]\in H_{b}^* \st f(e)=0\} 
\geq\frac{1}{2}
$$
for $j=0,\ldots ,l$.
Consider a continuous function $\chi : H^*_b \to [0,1]$
with compact support $K=\cl{\chi^{-1}(0,1]}$
disjoint from the closed subspace
$\{ [f,y]\in H_{b}^* \st f(e)=0\}$. 
It will be enough to show that $\int \chi\, {\rm d}\mu_j^b
\leq \frac{1}{2}$.
For $n$ sufficiently large, we have
$\{ [f,y]\in H^*_b \st f(e_n)\leq |y|/n\} \cap K=\emptyset$
and so $\int \chi\, {\rm d}\mu_j^{b_n}\leq \frac{1}{2}$.
But, by the continuity of the family,
$\int \chi\, {\rm d}\mu_j^{b_n}\to \int \chi\, {\rm d}\mu_j^b$
as $n\to\infty$.
\end{proof}
\end{appendix}

\end{document}